\documentclass[a4paper,11pt,reqno]{acta} %%based on AMS-LATEX

\usepackage{amssymb,latexsym}
\usepackage{amsmath,amsthm}
\usepackage{enumerate}
\usepackage[mathscr]{eucal}

\usepackage{graphics} %%EPS files of figures are accepted
\usepackage{graphicx}
\graphicspath{ {./images/} }

\theoremstyle{plain}

\theoremstyle{definition}

\theoremstyle{remark}

\allowdisplaybreaks

\begin{document}
\null
\vskip 1.8truecm

\title[Natural vibrations of curved nano-beams and nano-arches]
{Natural vibrations of curved nano-beams and nano-arches } 
\author{Jaan Lellep and Shahid Mubasshar}    %%%  by our standard, for multiple authors, the names of the authors must be written in alphabetical order %%%  
\address{Institute of Mathematics and Statistics, University of Tartu, Estonia}
\email{jaan.lellep@ut.ee} %%(if possible)}
%\urladdr{URL address (if need be)}
\address{Institute of Mathematics and Statistics, University of Tartu, Estonia}
\email{shahid.mubasshar@ut.ee}  %%(if possible)}
%\urladdr{URL address (if need be)}
%\address{etc.}

%\dedicatory{Dedication (if need be)}

\begin{abstract}
The natural vibrations of curved nano-beams and nano-arches are studied. The nano-arches under consideration have piece wise constant thickness; these are weakened with stable cracks located at re-entrant corners of the steps. A method of determination of natural frequencies is developed making use of the method of weightless rotating spring. The aim of the paper is to assess the sensitivity of the eigenfrequencies on the geometrical and physical parameters of the nano-arch. The results of the calculations favourably compare with similar works of other researchers.  
\end{abstract}

%\subjclass{2020 Mathematics Subject Classification number(s)}
%\keywords{Nano-arches, natural frequency, crack, step}
%\date{xxxx (will be added by the editors)}

%\thanks{https://doi.org/10.12697/ACUTM.xxxx.yy.zz (will be added later)}
%\thanks{Corresponding author: (if necessry)}  %% For thanks we suggest to use the section Acknowledgement(s) at the end of the paper  %%

\maketitle

%%%%%%%%%%%%%%%%%%%%%%%%%%%%%%%%%%%%%%%%%%%%%%%%%%%%%%%%%%%%%%%%%%%%%%%%%%%%%%%%%%%%%%%%%%%%%%%%%%%%%%%%%%%%%%%%%%%%%%%%%%%%%%%%%%%%%%%%%%%%%%%%%%%%%%%%%%%%%%%%%%%%%%%

\section{Introduction}
%%%%%%%%%%%%%%%%%%%%%%%%%%%%%%%%%%%%%%%%%%%%%%%%%%%%

The main ideas of the non-local theory of elasticity were initially formulated by Eringen \cite{7} and Eringen and Edelen \cite{8} several decades ago. However, the rapid progress in the non-local theories got its start with the wide use of nanomaterials (see Thai \cite{27}, Thai et al \cite{28}, Reddy \cite{21}). Comprehensive reviews of papers dedicated to the mechanical behaviour of nano-beams, nano-plates and nano-arches can be found in the review papers by Faghidian \cite{9,10} Farajpour, Ghoyesh and Farokhi \cite{11}, also in Wang and Arash \cite{31}. It is interesting to remark that the non-local theory was initially formulated in the integral form and later it was reformulated by Eringen in the differential form making use of a specific kernel function. Since the differential form is  more simple it is widely used in the analysis of nanostructures.

A general approach to the application of the non-local theories of the elasticity in the bending, buckling and vibration of nano-beams was developed by Aydogdu \cite{2}. Resorting to the results obtained by Nayfeh and Emam \cite{20} a comparative study on the buckling and post-buckling analysis of nano-beams was undertaken by Emam \cite{6}. The formulations corresponding to the classical Euler-Bernoulli approach, the first order Timoshenko theory and higher order shear deformation Reddy theories are analysed in detail. The comparison of classical bending theory and the shear deformation theories is carried out by Reddy \cite{21,22}, as well. 

Vibrations of nano-beams containing cracks and other defects are studied by Roostai and Haghpanahi \cite{24}, Loya \cite{19}, also by Lellep and Lenbaum\cite{17}, Loghmani and Yazdi \cite{18}, Hossain and Lellep \cite{12,13}. In the paper \cite{13} the effect of the temperature is taken into account. 

In the paper by Jiang and Wang \cite{14} analytical solutions are developed for nano-beams subjected to axial and thermal stresses. 

The vibrations of curved beams or arches and the detection of damages in these structures is studied by Cerri and Ruta \cite{3}, Viola, Hrtioli and Dilena \cite{29} with the help of vibrational analysis. Viola and Tournabene \cite{30} undertake the analysis of circular arches of variable cross-section. The effects of cracks are taken into account in the paper by Karaagae et al \cite{15}.

In the present study an attempt is made to define the natural frequencies of circular nano-beams and to study the sensitivity of eigenfrequencies on the geometrical and material parameters of the nano-arch.   
\section{Formulation of the problem}
Let us consider the dynamic behaviour of a nano-arch of radius $R$. For determination of the positions of central points of the arch the polar coordinates are introduced (Fig 1).
\begin{figure}[h!]
\centerline{\includegraphics[width=14cm]{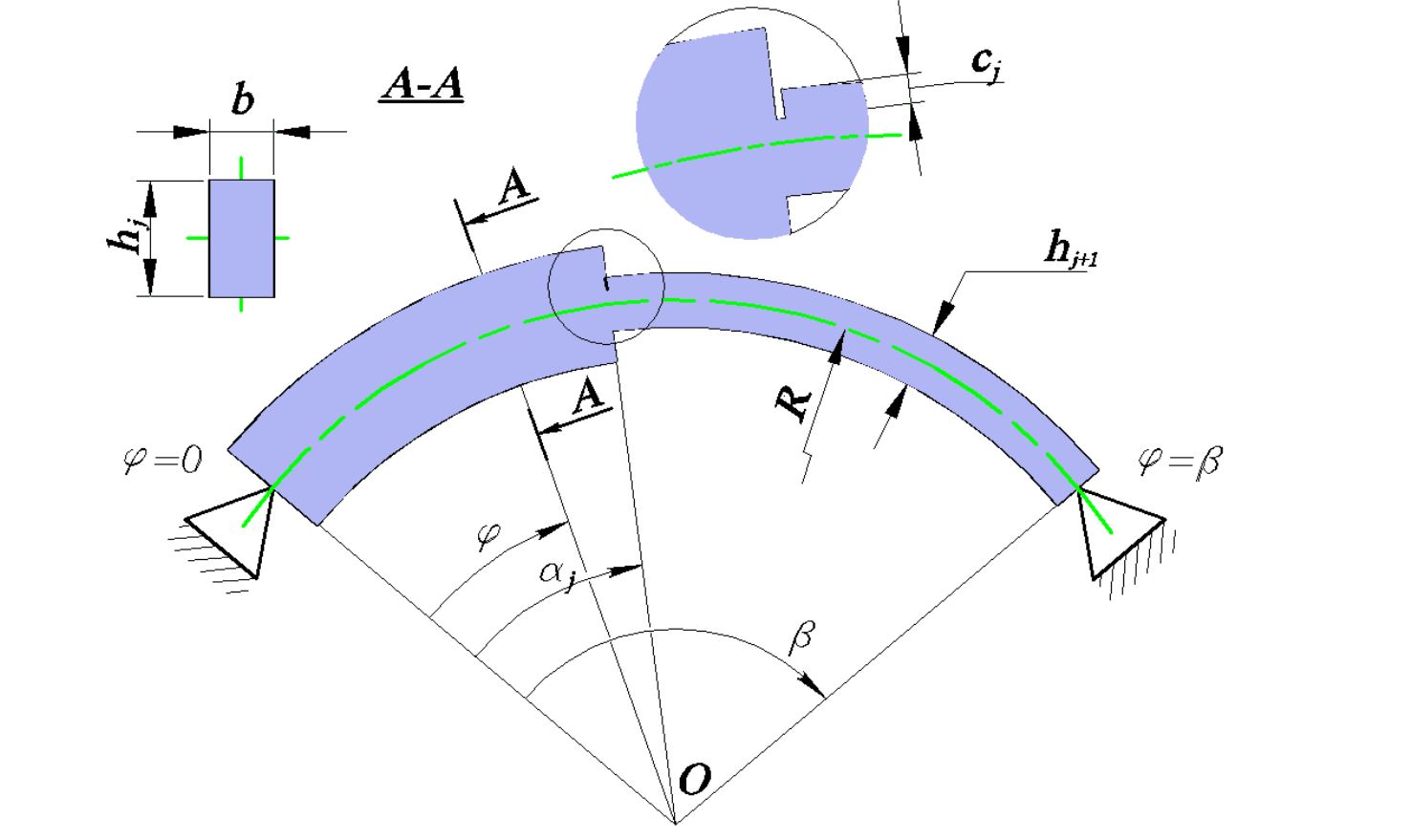}}
\caption{Elastic stepped nano-arch.}
\label{Fig.1}
\centering
\end{figure}
Assume that the angles $\varphi=0$ and $\varphi=\beta$ are correspond to the edges of the nano-arch, respectively. 

It is assumed that the thickness of the nano-arch is defined as 
\begin{equation}
\label{eqn:1}
h=
\begin{cases}
    \begin{alignedat}{4}
    &h_{0}\quad&,\quad& \varphi \in [0,\alpha)\\                        
    &h_{1}\quad&,\quad&\varphi \in (\alpha, \beta]\\
	\end{alignedat}
  \end{cases}
\end{equation}          
In (\ref{eqn:1}), $h_{0}, h_{1}$ are given numbers and $\alpha \in (0, \beta)$. The width $b$ of the nano-arch is assumed to be a constant. 
The aim of the study is to determine the eigenfrequencies of stepped nano-beams weakened with defects at the re-entrant corners of steps. It is assumed that at the cross section $\varphi=\alpha$ a crack with length $c$ is located. The defect is treated as a stable crack; no attention will be paid to the extension of the crack during vibrations. 

However, the sensitivity of the eigenfrequency with respect to the position of the crack and other parameters will be clarified. 
\section{Governing equations and main hypotheses}
The system of governing equations includes the equilibrium equations and constitutive and geometrical relations with boundary conditions. The equilibrium equations for an element of the arch or curved beam can be presented as (see Soedel \cite{25}),
\begin{equation}
\label{eqn:2}
\begin{aligned}
  \frac {\partial{N}} {\partial s}{}+{} \frac {Q}{R}&+p_{u}&= &&\bar \rho h\ddot{U},\\
 \frac {\partial{Q}} {\partial s} {}-{}\frac {N}{R}&+p&= &&\bar \rho h\ddot{W},\\
\frac {\partial{M}} {\partial{ s}}&-Q&=&&0
\end{aligned}
\end{equation}
In (\ref{eqn:2}), $N$, $M$ and $Q$ denote the membrane force, bending moment, and shear force, respectively. Here $s$ is curvilinear coordinate related to the current angle as
\begin{equation}
s=R\varphi.
\end{equation}
In  (\ref{eqn:2}), $U$, $W$ stand for the displacements in the circumferential and normal direction, respectively and $p_{u}$ and $p$ are the intensities of the distributed loading in these directions. Here $h$ stands for the thickness of the arch and $\bar \rho$ is the mass per unit length of the arch. Let $b$ be the width of the arch. In (\ref{eqn:2}), dots denote the differentiation with respect to time $t$. This means that
\begin{equation}
\ddot{U}= \frac {\partial^2{U}} {\partial t^2},\quad \ddot{W}=\frac {\partial^2{W}} {\partial t^2}.
\end{equation}
The strain-displacement equations can be taken according to Soedel \cite{25} as
\begin{equation}
\label{eqn:5}
\varepsilon=\frac {1} {R}W+\frac {\partial{U}} {\partial {s}}
\end{equation}
\begin{equation}
\label{eqn:6}
\varkappa=-\frac {\partial^2{U}} {\partial t^2}+\frac {\partial{U}} {\partial {s}.R}
\end{equation}
 In (\ref{eqn:5}),  (\ref{eqn:6}) $\varepsilon$ stands for the relative extension of the curved element and $\varkappa$ is the curvature of the middle line of the arch. 

It is assumed herein that the vibrations of the arch are in-place motions only. 

According to the Hook's law in the classical theory of the elasticity one has 
\begin{equation}
\label{eqn:7}
M_{c}=EI\varkappa,
\end{equation}
where $M_{c}$ denotes the bending moment in the classical theory Here $E$ is the Young modulus and $I$\textemdash the moment of inertia of the cross section of the arch. Evidently, in the case of the rectangular cross section 
\begin{equation}
\label{eqn:8}
I=\frac{1}{12}bh^3.
\end{equation}
As $h$ is a piecewise continuous function of the current angle the second moment (\ref{eqn:8}) has the values $I_{0}$ and $I_{1}$ in the regions ($0, \alpha$) and ($\alpha, \beta$), respectively.  
It is assumed that the material of the arch obeys the constitutive equations of the non-local theory of elasticity (see Eringen \cite{7}, Reddy \cite{21}). In the non-local theory the stresses at a given point of the body depend on the strains at all points of the body. One of the simplest physical relations of this type can be presented as (see Reddy \cite{21}, Emam \cite{6})
\begin{equation}
\label{eqn:9}
\sigma_{ij}-\eta\nabla^{2}\sigma_{ij}=\sigma^c_{ij},
\end{equation}
where $\sigma_{ij}$ denotes the components of non-local stress tensor and $\sigma^c_{ij}$ stands for the classical elastic stress tensor. Here $\eta=(e_{0}a)^{2}$ stands for a material constant, $a$ being the dimension of the lattice of the material and $e_{0}$\textemdash a physical constant to be defined experimentally.
\\
Instead of stress components $\sigma_{ij}$ in (\ref{eqn:9}) one can use the generalized stresses used in the equilibrium equations (\ref{eqn:2}).\\
However, it is reasonable to adapt the system (\ref{eqn:2}) to the current problem. Let us assume that $p_{u}=0$, $U(Q, t)$ is small and that $\ddot{U}$ is also negligible. In this case the system (\ref{eqn:2}) can be rewritten as 
\begin{equation}
\label{eqn:10}
\begin{aligned}
  N^{'}&=-Q,\\
  Q^{'}&=N+R(\rho h \ddot{W}-p),\\
 M^{'}&=RQ,
\end{aligned}
\end{equation}
where prim denotes the differentiation with respect to the current angle $\varphi$.\\
It can be seen from  (\ref{eqn:10}) that in the case when $M(\varphi )$ and $N(\varphi )$ vanish at the same cross section one has 
\begin{equation}
\label{eqn:11}
M=-RN.
\end{equation}
Equations  (\ref{eqn:10}) and  (\ref{eqn:11}) result in 
\begin{equation}
\label{eqn:12}
M^{''}+M+R^{2}(p-\rho h\ddot{W})=0.
\end{equation}
Taking in  (\ref{eqn:9}) instead of  $\sigma_{ij}$ the moment $M$ yields 
\begin{equation}
\label{eqn:13}
M-\eta M^{''}=M_{c}.
\end{equation}
Combining  (\ref{eqn:13}) and  (\ref{eqn:7}) with  (\ref{eqn:12}) one obtains
\begin{equation}
\label{eqn:14}
M=\frac{1}{1+\eta}[-EIW^{''}+h\eta \rho \ddot{W}]. 
\end{equation}
Substituting the bending moment from  (\ref{eqn:14}) into  (\ref{eqn:12}) one obtains the equation 
\begin{equation}
\label{eqn:15}
EI(W^{IV}+W^{''})+\eta \rho h R^{2} \ddot{W}-\eta \rho h R^{2}(1+\eta)\ddot{W}^{''}=0. 
\end{equation}
Equation (\ref{eqn:15}) is a fourth order equation. The boundary conditions for (\ref{eqn:15}) are in the case of simply supported arches
\begin{equation}
\label{eqn:16}
\begin{aligned}
  W(0, t)=0,\\
  W(\beta, t)=0
\end{aligned}
\end{equation}
and
\begin{equation}
\label{eqn:17}
\begin{aligned}
  M(0, t)=0,\\
  M(\beta, t)=0.
\end{aligned}
\end{equation}
\section{Solution of the governing equations}
Making use of the method of separation of variables it is assumed that 
 \begin{equation}
\label{eqn:18}
  W(\varphi, t)=X(\varphi).T(t),
\end{equation}
where $X(\varphi)$ a function of the coordinates $\varphi$ and $T(t)$\textemdash a function depending on time $t$ only. \\
Differentiating  (\ref{eqn:18}) one can find easily that
\begin{equation}
\label{eqn:19}
\begin{aligned}
  \ddot{W}&=X(\varphi).\ddot{T}(t),\\
  W^{''}&=X^{''}(\varphi).T(t),\\
W^{IV}&=X^{IV}(\varphi).T(t),\\
 \ddot{W}^{''}&=X^{''}(\varphi).\ddot{T}(t).
\end{aligned}
\end{equation}
The substitution of (\ref{eqn:19}) in (\ref{eqn:15}) leads to the differential equations
 \begin{equation}
\label{eqn:20}
  \ddot{T}+\omega ^{2}T=0
\end{equation}
and
 \begin{equation}
\label{eqn:21}
EI\frac{X^{IV}+X^{''}}{X}-\omega ^{2}\eta \rho h R^{2}(1-\frac{X^{''}}{X})=0
\end{equation}
In  (\ref{eqn:20}) $\omega$ stands for the frequency of natural vibrations of the arch. The equation  (\ref{eqn:21}) must be solved separately in regions $(0, \alpha)$ and $(\alpha, \beta)$ respectively, since $h$ is defined by  (\ref{eqn:1}) so that it has different values in these sections. Evidently, the equation  (\ref{eqn:20}) has periodical solutions. One of these can be taken as 
\begin{equation}
\label{eqn:22}
T=\bar A \sin (\omega t)
\end{equation}
Equation  (\ref{eqn:21}) can be presented in the form 
\begin{equation}
\label{eqn:23}
X^{''''}+(1+K_{j})X^{''}-K_{j}X=0,
\end{equation}
where
\begin{equation}
\label{eqn:24}
K_{j}=
\begin{cases}
    \begin{alignedat}{4}
    &K_{0}\quad&,\quad& \varphi \in (0,\alpha),\\                        
    &K_{1}\quad&,\quad&\varphi \in (\alpha, \beta).\\
	\end{alignedat}
  \end{cases}
\end{equation}
In  (\ref{eqn:23}),  (\ref{eqn:24})
\begin{equation}
\label{eqn:25}
K_{j}=\frac{\omega ^{2}\eta \rho h_{j} R^{2}}{EI_{j}}
\end{equation}
where $j=0, 1$. Here
\begin{equation}
\label{eqn:26}
I_{0}=\frac{bh_{0}^3}{12}\quad,\quad I_{1}=\frac{bh_{1}^3}{12}
\end{equation}
In order to solve the linear fourth order equation  (\ref{eqn:23}) one has to solve the characteristic equations 
\begin{equation}
\label{eqn:27}
\lambda_{j}^{4}+(1+K_{j})\lambda_{j}^{2}-K_{j}=0
\end{equation}
It is easy to recheck that 
\begin{equation}
\label{eqn:28}
\lambda_{j}^{2}=\frac{1}{2}(-1-K_{j}\pm\sqrt{1+6K_{j}+K_{j}^{2}})
\end{equation}
and 
\begin{equation}
\label{eqn:29}
\lambda_{j}=\pm\{\frac{-1-K_{j}}{2}\pm \frac{1}{2}\sqrt{1+6K_{j}+K_{j}^{2}}\}^{\frac{1}{2}}
\end{equation}
Thus the general solution of the (\ref{eqn:21}) has the form
\begin{equation}
\label{eqn:30}
X_{j}=C_{1} \cosh\mu_{j} \varphi+C_{2} \sinh\mu_{j} \varphi+C_{3} \cos\nu_{j} \varphi+C_{4} \sin\nu_{j}\varphi
\end{equation}
where
\begin{equation}
\label{eqn:31}
\begin{aligned}
 \mu_{j}&=-\frac{1}{2}(1+K_{j})+ \frac{1}{2}\sqrt{1+6K_{j}+K_{j}^{2}} \\
  \nu_{j}&=\frac{1+K_{j}}{2}+ \frac{1}{2}\sqrt{1+6K_{j}+K_{j}^{2}}
\end{aligned}
\end{equation}
 It is expected in  (\ref{eqn:27})- (\ref{eqn:31}) that $j=0$, if $\varphi\in (0, \alpha)$ and $j=1$ if $\varphi\in(\alpha, \beta)$.\\
The boundary conditions  (\ref{eqn:16}) and  (\ref{eqn:17}) together with  (\ref{eqn:18}) admit to assert that the boundary requirements for $X_{j}(t)$ are 
\begin{equation}
\label{eqn:32}
\begin{aligned}
 X_{j}(0)&=0, \\
 X^{''}_{j}(0)&=0,
\end{aligned}
\end{equation}
and
\begin{equation}
\label{eqn:33}
\begin{aligned}
 X_{j}(\beta)&=0, \\
 X^{''}_{j}(\beta)&=0.
\end{aligned}
\end{equation}
Note that (\ref{eqn:32}) and (\ref{eqn:33}) hold good in the case of the arch simply supported at both edges. 
\section{Local flexibility due to the crack}
It is recognized long time ago that cracks and other defects are the sources of additional structural compliance. Dimarogonas \cite{5} Chondros et al \cite{4} also Kukla \cite{16}, suggested to treat the slope of deflection $W^{'}$ as a discontinuous variable and to define a new variable
\begin{equation}
\label{eqn:34}
\theta=W^{'}(\alpha+, t)-W^{'}(\alpha-0, t)
\end{equation}
so that 
\begin{equation}
\label{eqn:35}
\theta=cM(\alpha, t)
\end{equation}
where $c$ can be considered as the additional compliance. In  (\ref{eqn:35}) $M$ stands for the moment of internal forces at $\varphi=\alpha$.
Evidently, $c$ can be an appropriate matrix, as well. In this case instead of $M$ one has the vector of corresponding internal forces applied  at the same cross section.\\
It is known in the linear elastic fracture mechanics that the energy release rate due to the crack can be calculated as 
\begin{equation}
\label{eqn:36}
G=\frac{M^{2}}{2b}\frac{dC}{c}.
\end{equation}
In  (\ref{eqn:36}) $b$ is the width of the crack of rectangular cross section and $M$ is the moment applied at this cross section and $c$ stands for the length of the crack.  \\
On the other hand, the stress intensity factor can be defined as (see Anderson \cite{1}),
\begin{equation}
\label{eqn:37}
K=\sigma \sqrt{\pi c} F(s)
\end{equation} 
where $s=\frac{c}{h}$ and 
\begin{equation}
\label{eqn:38}
\sigma=\frac{6M}{bh^{2}}
\end{equation}
where $F(s)$ is the shape function. \\
The quantity $K$ is coupled with the energy release rate as 
\begin{equation}
\label{eqn:39}
G=\frac{K^{2}}{E^{'}}
\end{equation}
In  (\ref{eqn:39}) $E^{'}=E$ for the plane stress state and  $E^{'}=\frac{E}{1-\nu^{2}}$ for plane strain state. Following Dimarogonas \cite{5} , Rizos et al \cite{23} one can use the function $F(s)$ in the form 
\begin{equation}
\label{eqn:40}
F(s)=1.93-3.07 s+14.53 s^{2}-25.11 s^{3}+25.8 s^{4}
\end{equation}
Another shape function often used by researchers is (see Tada et al \cite{25}, Dimarogonas \cite{5})
 \begin{equation}
\label{eqn:41}
F(s)=\frac{\sqrt {\tan\psi}}{\psi\cos\psi}[0.923+0.199(1-\sin \psi)^{4}]
\end{equation}
where $\psi=\frac{\pi s}{2}$.
\section{Numerical results}
Numerical results are presented for nano-arches simply supported at both ends in fig 2-6. 
Here the nano-arches with a single step whereas $R=30 nm$, $h_{0}=10 nm$, $h_{1}=20 nm$, the material constant are $E=7\times10^{11} Pa$, $\nu=0.3$
\begin{figure}[h!]
\centerline{\includegraphics[width=14cm]{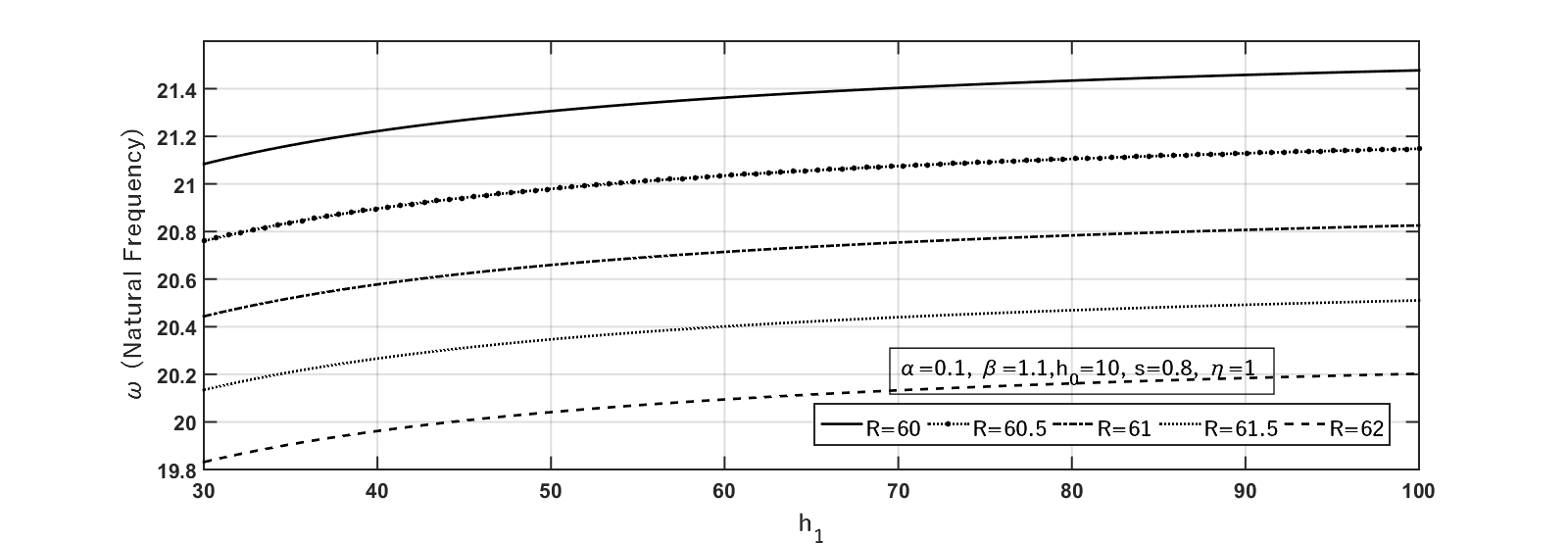}}
\caption{Natural frequency Vs thickness of the nano-arch}
\label{fig:fig.2}
\centering
\end{figure}
\begin{figure}[h!]
\centerline{\includegraphics[width=14cm]{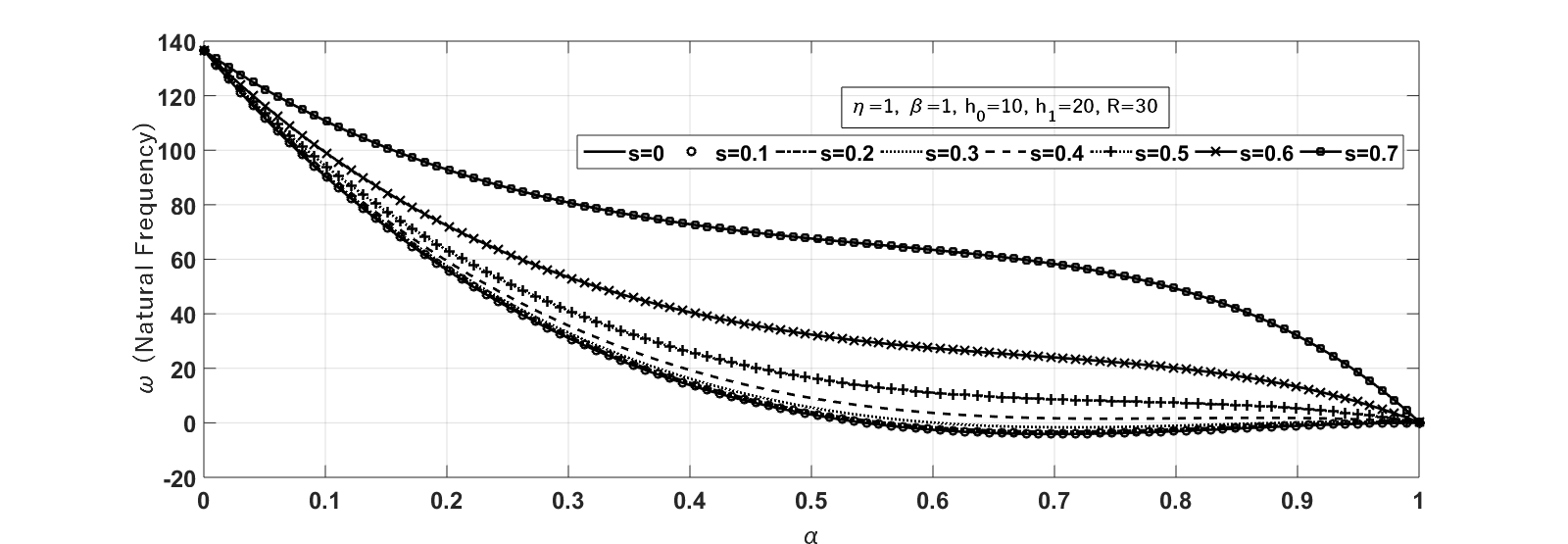}}
\caption{Natural frequency Vs step location}
\label{fig:fig.3}
\centering
\end{figure}
\begin{figure}[h!]
\centerline{\includegraphics[width=14cm]{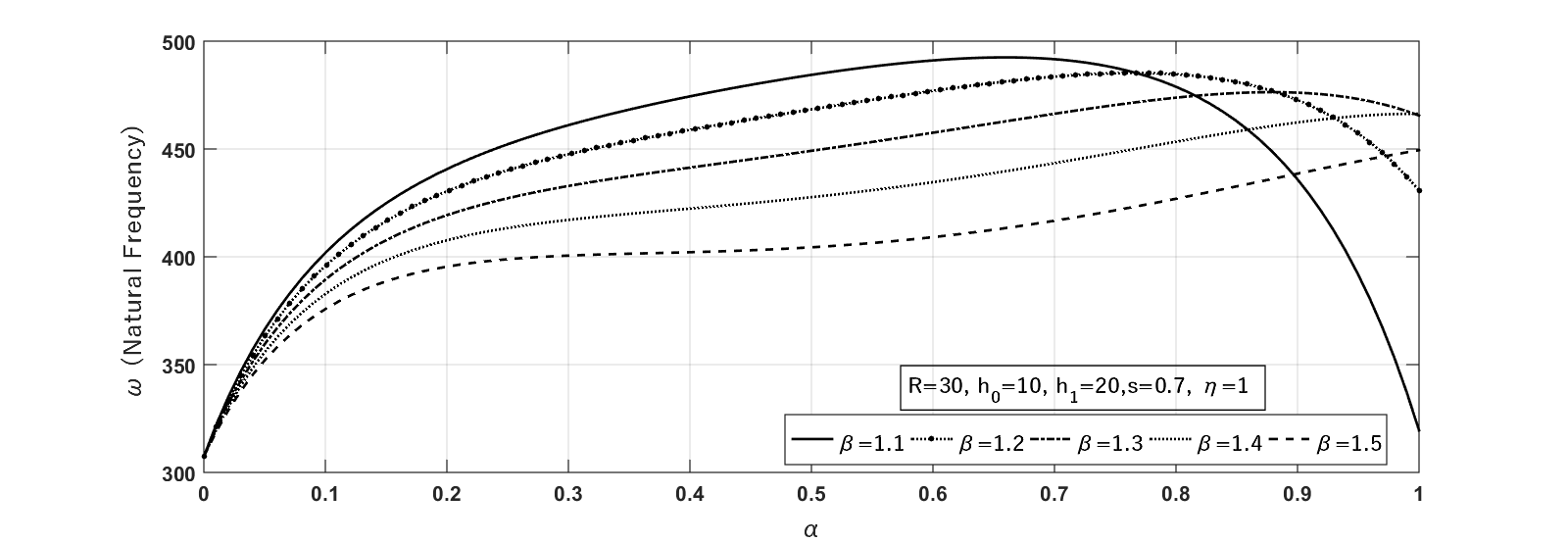}}
\caption{Natural frequency Vs step location for different values of $\beta$}
\label{fig:fig.4}
\centering
\end{figure}
\begin{figure}[h!]
\centerline{\includegraphics[width=14cm]{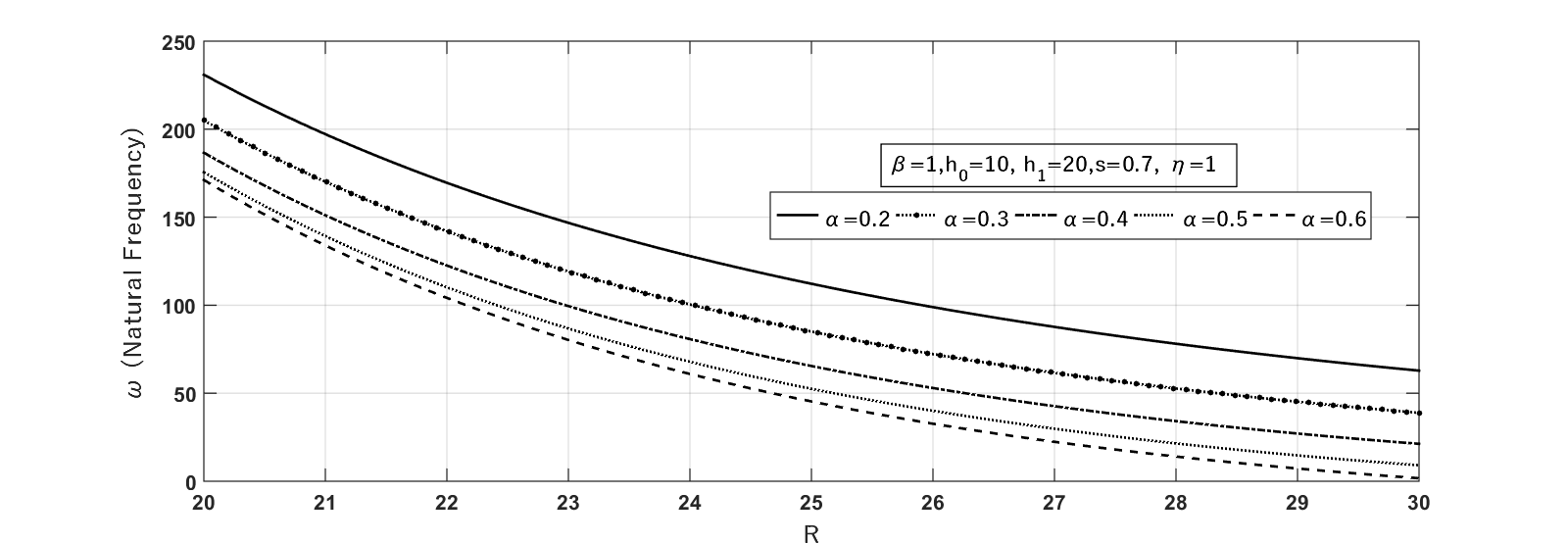}}
\caption{Natural frequency Vs radius of the nano-arch.}
\label{fig:fig.5}
\centering
\end{figure}
\begin{figure}[h!]
\centerline{\includegraphics[width=14cm]{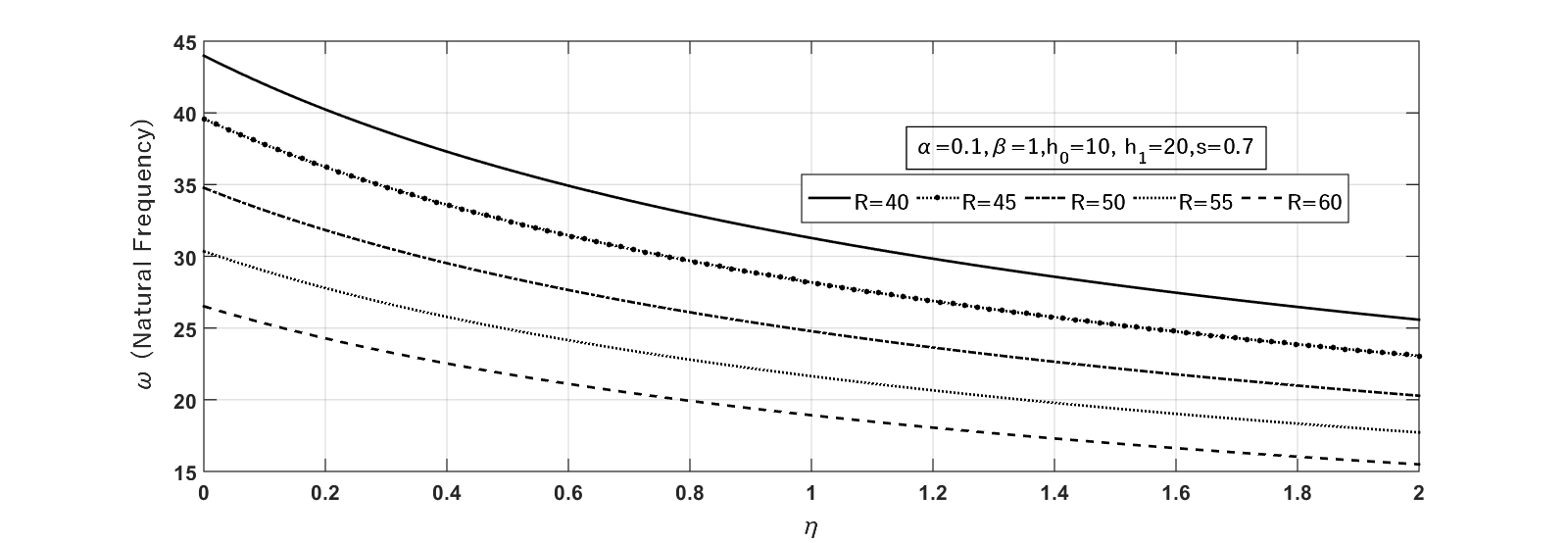}}
\caption{Natural frequency Vs material constant of the nano-arch.}
\label{fig:fig.6}
\centering
\end{figure}
\begin{figure}[h!]
\centerline{\includegraphics[width=14cm]{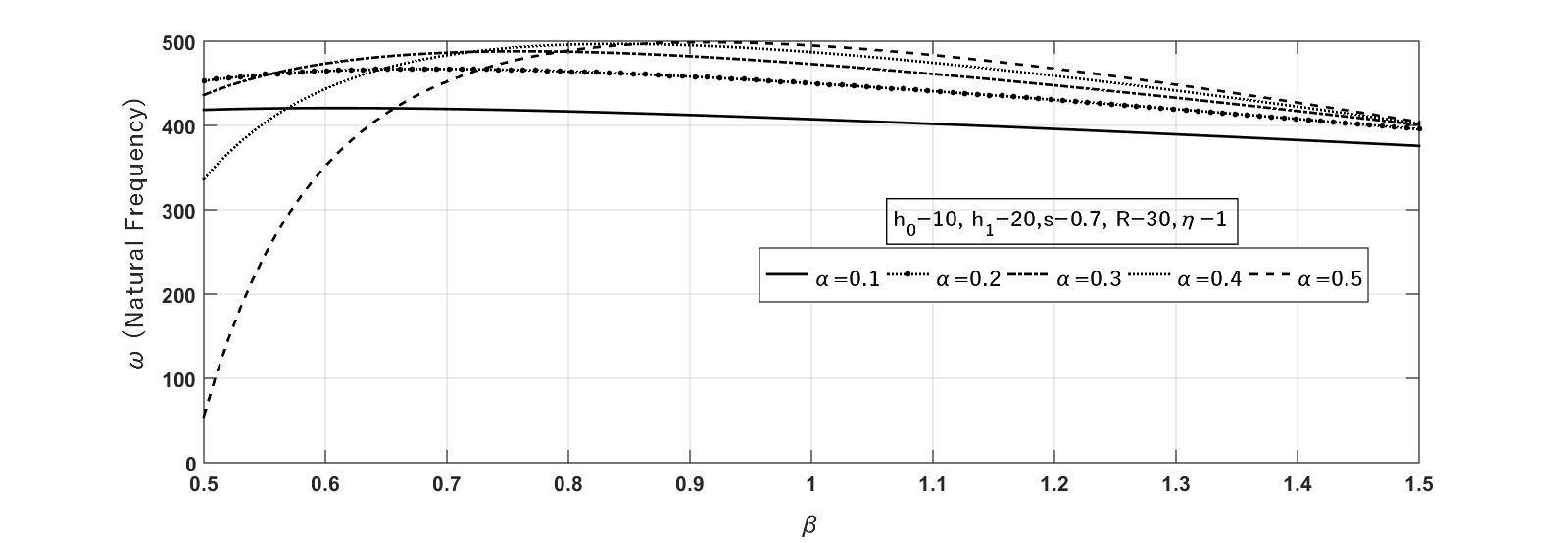}}
\caption{Natural frequency Vs central angle of the nano-arch.}
\label{fig:fig.7}
\centering
\end{figure}
\begin{table}[t]
\centering
\begin{tabular}{ |c|c|c|c| } 
 \hline
 Mode&$\eta$ & Present & Thai\cite{27} \\ 
  \hline
 1&0& 9.75821445 & 9.2745 \\ 
  
& 1 & 7.05584295 & 8.8482 \\ 
 
 &2& 5.80188130 & 8.4757 \\ 
 
 &3 & 5.04192108 & 8.1466 \\ 
 
 &4 & 4.51883759 & 7.8530 \\ 
 \hline
\end{tabular}

\caption{Natural vibration of the nano-arches of the constant thickness, without any crack and with central angle very very small}
\label{tab:1}
\end{table}

In figure (\ref{fig:fig.2}) the natural frequency of the nano-arch versus the thickness is depicted for different values of the radius $R$. It can be seen from the figure (\ref{fig:fig.2}) that the natural frequency increases if the thickness of the arch increases. However, if the radius of the arch increases then the natural frequency decreases.\\
In the next figures (\ref{fig:fig.3})-(\ref{fig:fig.7}) corresponding data for one-stepped arches is presented. Here $h_{0}=10 nm$,  $h_{1}=20 nm$.\\
In figures (\ref{fig:fig.3}) and (\ref{fig:fig.4}) the natural frequency as a function of the step location $\alpha$ is presented. In figure (\ref{fig:fig.3}) the nano-arch with the central angle $\beta=1 rad$ is treated. Different curves in figure (\ref{fig:fig.3}) are correspond to the crack extensions $s=0$, $s=0.1$. $s=0.2$, $s=0.3$, $s=0.4$,$s=0.5$,$s=0.6$,$s=0.7$ respectively. It can be seen from figure (\ref{fig:fig.3}) that the lowest values of the natural frequency correspond to the arch without any cracks. 
In figure (\ref{fig:fig.4}) similar results are presented for different angle values of the central angle $\beta$ (here $s=0.7$).
It reveals from figure (\ref{fig:fig.4}) that in the case of smaller values of the angle $\alpha$ the lowest eigenfrequency is achieved in the case of the largest value of the central angle $\beta$. However, in the case when $\beta >2\alpha$ this relationship is more complicated.\\
The relationship between the natural frequency and the radius $R$ of the nano-arch is shown in figure (\ref{fig:fig.5}) in the case $\beta=1$ and $s=0.7$ (here $\eta=1 nm$). \\
It reveals from figure (\ref{fig:fig.5}) that the smaller is $\alpha$ the larger is the the natural frequency of the nano-arch.\\
On the other hand, the larger is $\alpha$ the lower is the eigenfrequency. \\
The dependence of the eigenfrequency on the material parameter $\eta=(e_{0}a)^{2}$ is demonstrated in figure (\ref{fig:fig.6}) for different value of the radius $R$. Generally speaking, the eigenfrequency decreases if the material parameter increases. Thus, in the case $\eta=0$ the natural frequency has its maximal value.\\
The sensitivity of the eigenfrequency on the central angle $\beta$ of the nano-arch is shown in figure (\ref{fig:fig.7}) for radius $R=30nm$. Different curves in figure (\ref{fig:fig.7}) correspond to nano-arches having the step in different cross sections (here $\eta=1 nm$). One can see from figure (\ref{fig:fig.7}) that in the case of larger values of the angle $\beta$ curves corresponding to different values of the step angle are quite close to each other. However, if $\beta$ is around the value $0.5$ the discrepancies between these curves are large. \\
The results obtained in the current study are compared with those obtained by Thai \cite{27} in the case when $\alpha=0$, $\beta=0.5^{\circ}$. Table (1) shows that the results are in reasonable correspondence. 
\section{Concluding remarks}
Natural vibrations of curved nano-beams are treated in the current study. The nano-beams have piecewise constant thickness with stable cracks at the re-entrant corners of steps. A method for determination of natural frequencies based on the methods of the linear elastic fracture mechanics is developed. Calculations carried out showed that the cracks effect essentially the values of natural frequencies. \\
It is shown that the lowest values of the natural frequency correspond to nano-arches without any defect. The results of current study are compared with results obtained by Thai\cite{27}. The results compare favourably in the case of small values of the parameter $\eta$.  

\bibliographystyle{vancouver}
\bibliography{ACUTM}

\begin{thebibliography}{10}

\bibitem{7}
Eringen AC, Wegner J.
\newblock Nonlocal continuum field theories.
\newblock Appl Mech Rev. 2003;56(2):B20--B22.

\bibitem{8}
Eringen AC, Edelen D.
\newblock On nonlocal elasticity.
\newblock International journal of engineering science. 1972;10(3):233--248.

\bibitem{27}
Thai HT.
\newblock A nonlocal beam theory for bending, buckling, and vibration of
  nanobeams.
\newblock International Journal of Engineering Science. 2012;52:56--64.
\newblock Available from:
  \url{https://www.sciencedirect.com/science/article/pii/S0020722511002412}.

\bibitem{28}
Thai HT, Vo TP, Nguyen TK, Kim SE.
\newblock A review of continuum mechanics models for size-dependent analysis of
  beams and plates.
\newblock Composite Structures. 2017;177:196--219.

\bibitem{21}
Reddy J.
\newblock Nonlocal theories for bending, buckling and vibration of beams.
\newblock International journal of engineering science. 2007;45(2-8):288--307.

\bibitem{9}
Faghidian SA.
\newblock On non-linear flexure of beams based on non-local elasticity theory.
\newblock International Journal of Engineering Science. 2018;124:49--63.

\bibitem{10}
Ali~Faghidian S.
\newblock Unified formulations of the shear coefficients in Timoshenko beam
  theory.
\newblock Journal of Engineering Mechanics. 2017;143(9):06017013.

\bibitem{11}
Farajpour A, Ghayesh MH, Farokhi H.
\newblock A review on the mechanics of nanostructures.
\newblock International Journal of Engineering Science. 2018;133:231--263.

\bibitem{31}
Wang Q, Arash B.
\newblock A review on applications of carbon nanotubes and graphenes as
  nano-resonator sensors.
\newblock Computational Materials Science. 2014;82:350--360.

\bibitem{2}
Aydogdu M.
\newblock A general nonlocal beam theory: its application to nanobeam bending,
  buckling and vibration.
\newblock Physica E: Low-dimensional Systems and Nanostructures.
  2009;{\bf41}(9):1651--1655.

\bibitem{20}
Nayfeh AH, Emam SA.
\newblock Exact solution and stability of postbuckling configurations of beams.
\newblock Nonlinear Dynamics. 2008;54(4):395--408.

\bibitem{6}
Emam SA.
\newblock A general nonlocal nonlinear model for buckling of nanobeams.
\newblock Applied Mathematical Modelling. 2013;37(10-11):6929--6939.

\bibitem{22}
Reddy J.
\newblock Nonlocal nonlinear formulations for bending of classical and shear
  deformation theories of beams and plates.
\newblock International Journal of Engineering Science. 2010;48(11):1507--1518.

\bibitem{24}
Roostai H, Haghpanahi M.
\newblock Vibration of nanobeams of different boundary conditions with multiple
  cracks based on nonlocal elasticity theory.
\newblock Applied Mathematical Modelling. 2014;38(3):1159--1169.

\bibitem{19}
Loya J, L{\'o}pez-Puente J, Zaera R, Fern{\'a}ndez-S{\'a}ez J.
\newblock Free transverse vibrations of cracked nanobeams using a nonlocal
  elasticity model.
\newblock Journal of applied physics. 2009;105(4):044309.

\bibitem{17}
Lellep J, Lenbaum A.
\newblock Natural vibrations of stepped nanobeams with defects.
\newblock Acta et Commentationes Universitatis Tartuensis de Mathematica.
  2019;23(1):143--158.

\bibitem{18}
Loghmani M, Yazdi MRH.
\newblock An analytical method for free vibration of multi cracked and stepped
  nonlocal nanobeams based on wave approach.
\newblock Results in Physics. 2018;11:166--181.

\bibitem{12}
Hossain M, Lellep J.
\newblock Natural vibrations of nanostrips with cracks.
\newblock Acta et Commentationes Universitatis Tartuensis de Mathematica.
  2021;(1):87--105.

\bibitem{13}
Hossain M, Lellep J.
\newblock Natural vibration of stepped nanoplate with crack on an elastic
  foundation.
\newblock In: IOP Conference Series: Materials Science and Engineering. vol.
  660. IOP Publishing; 2019. p. 012051.

\bibitem{14}
Jiang J, Wang L.
\newblock Analytical solutions for thermal vibration of nanobeams with elastic
  boundary conditions.
\newblock Acta mechanica solida sinica. 2017;30(5):474--483.

\bibitem{3}
Cerri MN, Ruta GC.
\newblock Detection of localised damage in plane circular arches by frequency
  data.
\newblock Journal of Sound and Vibration. 2004;270(1-2):39--59.

\bibitem{29}
Viola E, Artioli E, Dilena M.
\newblock Analytical and differential quadrature results for vibration analysis
  of damaged circular arches.
\newblock Journal of Sound and Vibration. 2005;288(4-5):887--906.

\bibitem{30}
Viola E, Tornabene F.
\newblock Vibration analysis of damaged circular arches with varying
  cross-section.
\newblock Structural Durability \& Health Monitoring. 2005;1(2):155.

\bibitem{15}
Karaagac C, Ozturk H, Sabuncu M.
\newblock Crack effects on the in-plane static and dynamic stabilities of a
  curved beam with an edge crack.
\newblock Journal of Sound and Vibration. 2011;330(8):1718--1736.

\bibitem{25}
Soedel W.
\newblock Vibrations of shells and plates.
\newblock CRC Press; 2004.

\bibitem{5}
Dimarogonas AD.
\newblock Vibration of cracked structures: a state of the art review.
\newblock Engineering fracture mechanics. 1996;55(5):831--857.

\bibitem{4}
Chondros T, Dimarogonas A, Yao J.
\newblock A continuous cracked beam vibration theory.
\newblock Journal of sound and vibration. 1998;215(1):17--34.

\bibitem{16}
Kukla S.
\newblock Free vibrations and stability of stepped columns with cracks.
\newblock Journal of Sound and Vibration. 2009;319(3-5):1301--1311.

\bibitem{1}
Anderson TL.
\newblock Fracture mechanics: fundamentals and applications.
\newblock CRC press; 2017.

\bibitem{23}
Rizos P, Aspragathos N, Dimarogonas A.
\newblock Identification of crack location and magnitude in a cantilever beam
  from the vibration modes.
\newblock Journal of sound and vibration. 1990;138(3):381--388.

\end{thebibliography}

\end{document}